\documentclass{article}

\usepackage{amsmath,amssymb,amscd,theorem,pb-diagram}

{\theorembodyfont{\slshape}\newtheorem{theorem}{Theorem}}
{\theorembodyfont{\slshape}\newtheorem{corollary}[theorem]{Corollary}}
{\theorembodyfont{\slshape}\newtheorem{lemma}[theorem]{Lemma}}
{\theorembodyfont{\slshape}\newtheorem{note}[theorem]{Note}}

\newcommand{\Fp}{\ensuremath{\mathbb{F}_p}}
\renewcommand{\O}{\ensuremath{\mathcal{O}}}
\newcommand{\Ot}{\ensuremath{\widetilde{\mathcal{O}}}}
\newcommand{\G}{\ensuremath{\Gamma}}
\newcommand{\g}{\ensuremath{\gamma}}
\newcommand{\og}{\ensuremath{\bar{\gamma}}}
\newcommand{\Q}{\ensuremath{\mathbb{Q}}}
\newcommand{\Qp}{\ensuremath{\mathbb{Q}_p}}

\newcommand{\Zp}{\ensuremath{\mathbb{Z}_p}}

\newcommand{\R}{\ensuremath{\mathbb{R}}}
\newcommand{\Z}{\ensuremath{\mathbb{Z}}}

\newcommand{\F}{\ensuremath{\mathbb{F}}}

\newcommand{\ord}{\ensuremath{\textnormal{ord}}}

\newcommand{\K}{\ensuremath{\mathcal{K}}}
\newcommand{\C}{\ensuremath{\mathcal{C}}}
\newcommand{\D}{\ensuremath{\mathcal{D}}}
\renewcommand{\L}{\ensuremath{\mathcal{L}}}
\newcommand{\ep}{\hfill $\blacksquare$\vspace{\baselineskip}} 

\begin{document}

\title{Memory efficient hyperelliptic curve point counting}

\author{ Hendrik Hubrechts\\
\small{Research Assistant of the Research Foundation - Flanders (FWO - Vlaanderen)}\\
\small{Department of mathematics, Katholieke Universiteit Leuven}\\
\small{Celestijnenlaan 200B - bus 2400, 3001 Leuven (Belgium)}\\
\small{\texttt{Hendrik.Hubrechts@wis.kuleuven.be}}}
\date{February 19, 2010}

\maketitle

\begin{abstract}
In recent algorithms that use deformation in order to compute the number of points on varieties over a finite field, certain differential equations of matrices over $p$-adic fields emerge. We present a novel strategy to solve this kind of equations in a memory efficient way. The main application is an algorithm requiring quasi-cubic time and only quadratic memory in the parameter $n$, that solves the following problem: for $E$ a hyperelliptic curve of genus $g$ over a finite field of extension degree $n$ and small characteristic, compute its zeta function. This improves substantially upon Kedlaya's result which has the same quasi-cubic time asymptotic, but requires also cubic memory size.
\end{abstract}

\vspace{\baselineskip}

\section{Introduction and results}\label{sec:intro}
Originally motivated by cryptography (see \cite{CohenFrey} for an overview), in recent years much effort was put in finding algorithms that compute zeta functions of varieties over finite fields. The most efficient algorithms often use deformation, i.e.\ they `deform' the input variety to another variety that is easier to handle. This use of deformation originates from the work of Lauder \cite{LauderDeformation} in computing zeta functions of higher dimensional varieties. In that paper, the deformation step allowed him to reduce the dependency on the dimension in the algorithms. Also Tsuzuki \cite{TsuzukiKloosterman} came up with this idea in the context of the computation of Kloosterman sums. Later Gerkmann \cite{GerkmannEC} and the present author \cite{HubrechtsHECOdd,HubrechtsHECEven} showed that even in dimension one, profit can be drawn from the use of deformation.

Central in all these applications stands a certain $p$-adic differential equation, namely the Picard-Fuchs equation of the associated connection. In the present paper we give an algorithm that allows us to compute a particular solution of this equation in a memory efficient way. Combining this with a well chosen deformation, from a general hyperelliptic curve to one defined over the prime field, yields our main result (Theorem \ref{thm:OddChar}), stating that for a hyperelliptic curve of genus $g$ over the finite field $\F_{p^n}$, $p$ odd, the zeta function can be computed in time and memory (where we assume $p$ fixed and count bit operations)
\[\Ot(n^3g^{6.376})\quad\text{respectively}\quad \O(n^2g^4(\log g)^{2}).\]
This result can be compared to Kedlaya's algorithm \cite{KedlayaCountingPoints} that has $\Ot(n^3g^4)$ as time and $\Ot(n^3g^3)$ as space requirements. The crucial improvement is hence that our algorithm only requires an amount of memory quadratic in $n$. Note that the $\Ot$ is the soft-Oh notation as defined in \cite[Definition 25.8]{ModernCompAlg} (it is essentially a big-Oh notation that ignores logarithmic factors).

Later on Denef and Vercauteren \cite{DenefVercauteren} extended Kedlaya's algorithm to characteristic two. Combining this with our new result yields `on average' the same result as in odd characteristic, see Theorem \ref{thm:EvenChar} in Section \ref{sec:evenChar}. The result for the `general case' is slightly worse.

For any small characteristic it is also possible to compute the zeta functions of $n$ curves within a one dimensional linear family all at once in time $\O(n^{3+\rho})$ for arbitrary $\rho>0$. This result is presented in Section \ref{ssec:manyCurves}, and in Section \ref{ssec:hypersurfaces} an additional application concerning hypersurfaces is explained. We note that all our complexities are bitwise unless mentioned otherwise.

Before we prove these results, we give in Section \ref{sec:diffEq} a very general form of the algorithms involved, including a thorough investigation of the error propagation during the computations, all of which is concluded in Theorem \ref{thm:solvingDiffEq}.

The author wishes to thank Wouter Castryck, Filip Cools, Jan Denef, Jan Tuitman and the referee for their helpful remarks.

\section{The differential equation}\label{sec:diffEq}
In this section we will define the differential equation referred to above. Given some conditions on the coefficients and on certain local solutions of this equation, we can present two algorithms that solve it, together with their complexity analysis.

\subsection{A general kind of $p$-adic differential equation}\label{ssec:diffEq}

Let $p$ be a prime number and $\mathbb{K}$ a degree $n$ field extension of the field of $p$-adic numbers $\Q_p$. We denote with $\ord$ the valuation on $\mathbb{K}$ normalized to $\ord(p)=1$, and $\mathcal{O}_{\mathbb{K}}:=\{x\in \mathbb{K}\,|\,\ord(x)\geq 0\}$ is the ring of integers of $\mathbb{K}$. Let $m>0$ be an integer. If we say that we are working in $\mathbb{K}$ modulo $p^m$, we mean that we use absolute precision, i.e.\ two numbers are considered equal when their difference has valuation at least $m$. We will also use power series in $\G$ modulo $\G^\ell$ for some integer $\ell>0$. The dimension of the square matrices that we will encounter is denoted by $d$. If $A(\G)$ is a ($d\times d$) matrix over $\mathbb{K}[[\G]]$, we will always use $A_i\in\mathbb{K}^{d\times d}$ in the following sense: $A(\G)=\sum_iA_i\G^i$, and in order to ease notation we will often write $A$ instead of $A(\G)$. The valuation $\ord(A(\G))=\ord(A)$ is defined to be $\inf_i(\ord(A_i))$ when this infimum exists, and $-\infty$ otherwise. We say that a series $\sum_iA_i\G^i$ with $A_i\in \mathbb{K}^{d\times d}$ is $(x,y)$-log convergent for real numbers $x\geq 0$ and $y$ if for every $i\geq 0$
\[\ord(A_i)\geq -x\lceil\log_p(i+1)\rceil-y.\] This implies in particular that such a series converges on the open unit disk in $\mathbb{K}$.
The following easy lemma can be found as Lemma 15 in \cite{HubrechtsHECEven}.
\begin{lemma}\label{lem:convergenceLogProduct}
If $\sum_i A_i\G^i$ and $\sum_iB_i\G^i$ converge $(x,y)$-log resp.\ $(x',y')$-log, then their product has $(x+x',y+y')$-log convergence.\end{lemma}

Let $A, B, X, Y$ be matrices over $\mathbb{K}[\G]$ such that $A_0$ and $B_0$ are invertible. We define $\Delta$ as the $\mathbb{K}$-linear operator acting on $\mathbb{K}[[\G]]^{d\times d}$ by
\begin{equation}\label{eq:Delta}
K\quad\mapsto\quad\Delta K=A \frac{d K}{d \G} B+A K X+Y K B.
\end{equation}

\subsection{Requirements for the equation}\label{ssec:assumptions}

In the proof of Theorem \ref{thm:solvingDiffEq} below, we show that for every boundary condition $K_0$, a unique solution $K\in\mathbb{K}[[\G]]^{d\times d}$ of $\Delta K=0$ exists. Our goal is to compute an approximation modulo $p^m$ of this unique solution $K$, where we assume that we know $K_0,A,B,X$ and $Y$ up to arbitrary precision, and in addition that $K_0$ (and hence $K(\G)$) are invertible. Write $\zeta$ for an upper bound for the following three values:
\[\deg A+\deg B,\ \ \deg A+\deg X+1,\ \ \deg Y+\deg B+1.\]
We note that $A(\G)$ and $B(\G)$ are invertible in $\mathbb{K}[[\G]]^{d\times d}$ and require that there exists some  $\alpha\in\R_{\geq 0}$ such that
\[\ord(A),\ \ord(A^{-1}),\ \ord(B),\ \ord(B^{-1}),\ \ord(K),\ \ord( K^{-1})\geq -\alpha.\]
We write $C$ for the unique solution of $A\frac{d C}{d\G}+YC=0$, $C_0=1$, over $\mathbb{K}[[\G]]^{d\times d}$ and $D$ for the solution of $\frac{dD}{d\G}B+DX=0$, $D_0=1$. Then we require the existence of constants $\g\geq 0$ and $\delta$ such that
\begin{equation}\label{eq:CDlogConv}
C,\ C^{-1},\ D,\ D^{-1}\text{ have }(\g,\delta)\text{-log convergence.}\end{equation}
Finally we define $\psi:=5(\alpha+\delta)$ and we require that $\ord(X),\ \ord(Y)\geq -\psi$. As a last assumption we need that $K(\G)$ modulo $p^m$ consists of polynomials of degree less than $\ell$, so that our approximation will be a finite object.

Note that although these conditions combined seem to be quite severe, they are met by a suitably adapted version of the differential equations appearing in the intended point counting algorithms.

\subsection{Solving the differential equation}\label{ssec:solution}

Due to precision loss during the execution of the algorithm, we have to work initially with a higher precision than $p^m$. Therefore we define $\varepsilon := m+(5\g+1)\lceil\log_p\ell\rceil+\psi$ and we work modulo $p^\varepsilon$. Let $\omega$ be an exponent for matrix multiplication, meaning that we can compute the product of two $d\times d$ matrices over $\mathbb{K}$ using $\O(d^\omega)$ arithmetic operations in $\mathbb{K}$. We may take $\omega = 2.376$, see \cite{CoppersmithWinograd}.
\begin{theorem}\label{thm:solvingDiffEq}
Suppose that we know $K_0,A,A_0^{-1},B,B_0^{-1},X$ and $Y$ modulo $p^\epsilon$ and that all assumptions of Section \ref{ssec:assumptions} are met. Then we can compute $K(\G)$ modulo $p^m$ using $\Ot(\ell\zeta d^\omega n \varepsilon\log_2p)$ bit operations and with memory requirements $\O(\ell d^2n\varepsilon\log_2p)$ bits. Moreover, we can compute $K(1)$ modulo $p^m$ in the same amount of time but with memory only $\O(\zeta d^2n\varepsilon\log_2p)$ bits.\end{theorem}
\textsc{Proof.} We first give the algorithm, then we will determine how much precision is lost throughout the computations, and finally we will do a resource analysis.\\

Let $\K(\Gamma)$ and $\K(1)$ denote the approximations of $K(\G)$ resp.~$\K(1)$ that we compute using the following algorithm. Define the operator
$\Delta'$ on $\mathbb{K}[[\G]]^{d\times d}$ by
\begin{equation*}\label{eq:DeltaPrime}
K\quad\mapsto\quad\Delta' K=A_0^{-1}A \frac{d K}{d \G} BB_0^{-1}+A_0^{-1}A K XB_0^{-1}+A_0^{-1}Y K BB_0^{-1},
\end{equation*}
then clearly the equation $\Delta K=0$ is equivalent to $\Delta'K=0$. Writing down the coefficient of $\G^k$ in $\Delta' K=0$ gives an equality of the form
\[\sum_{a+b+c=k}\left( (b+1)A_0^{-1}A_aK_{b+1}B_cB_0^{-1}+A_0^{-1}A_aK_bX_cB_0^{-1}+A_0^{-1}Y_aK_bB_cB_0^{-1} \right)=0.\]
If we isolate $K_{k+1}$, this yields
\begin{equation}\label{eq:inductiveSoln}
(k+1)K_{k+1} = f_k\left(K_k,K_{k-1},\ldots,K_{k-(\zeta-1)}\right)
\end{equation}
for some easy to construct linear polynomial $f_k$ defined over $\mathbb{K}^{d\times d}$ (where we put $K_i=0$ for $i<0$). This recursion relation allows us to calculate $\K(\G)$ as follows: put $\K_i:=0$ for $i<0$, $\K_0:=K_0 \bmod p^\varepsilon$ and for $k=0,1,\ldots,\ell-2$ compute
\begin{equation}\label{eq:inductiveComp}
\K_{k+1} := \left[\frac 1{k+1}f_k\left(\K_k,\K_{k-1},\ldots,\K_{k-(\zeta-1)}\right)\right]\bmod p^\varepsilon.
\end{equation}
Finally we define $\K(\G):=\sum_{i=0}^{\ell- 1}\K_i\G^i$. Note that equation (\ref{eq:inductiveSoln}) implies that the solution $K(\G)$ of $\Delta K=0$ exists and is unique.

Computing $\K(1)$ uses the same idea, but in order to save memory we work as follows. Define $\L_0:=\K_0$ at the start, and put $\L_k := \L_{k-1}+\K_k\bmod p^{\varepsilon}$ each time that a new $\K_k$ is computed. This way we only have to store the last $\zeta$ matrices $\K_k,\K_{k-1},\ldots,\K_{k-(\zeta-1)}$ and $\L_k$ in order to find $\K_{k+1}$ and $\L_{k+1}$. After $\ell-1$ steps we will end with $\L_{\ell-1}=\sum_{i=0}^{\ell-1}\K_i\mod p^\varepsilon$.\\

Every time that we multiply non integral elements of $\mathbb{K}$, we can expect a certain loss in $p$-adic precision. We will now show that working modulo $p^\varepsilon$ suffices to conclude that $\K(\G)\equiv K(\G)\bmod p^m$ and $\K(1)\equiv K(1)\bmod p^m$. We follow the reasoning of Section 3.5 in \cite{HubrechtsHECEven}, but now in a more general setting\footnote{Gerkmann \cite{GerkmannHypersurfaces} independently has found a similar control of the error propagation.}.

With the appropriate recursion modulo $(p^\varepsilon,\G^\ell)$ similar to (\ref{eq:inductiveComp}), but now for $A_0^{-1}A\frac{dC}{d\G}+A_0^{-1}YC=0$ and $C_0=1$ we can compute $\C$, and in an analogous way $\D$ as approximation to $D$. We remark that $\mathcal{C}$ and $\mathcal{D}$ are only needed for the analysis, not in the actual computation. Rewriting (\ref{eq:inductiveComp}) implies that
\[
(k+1)\K_{k+1} - f_k\left(\mathcal{K}_k,\mathcal{K}_{k-1},\ldots, \mathcal{K}_{k-(\zeta-1)}\right) = p^\varepsilon\cdot(\text{integral error matrix}),\]
and if we sum over all $k$ this gives $\Delta' \mathcal{K}=p^\varepsilon\mathcal{E}_\K$, where $\mathcal{E}_\K$ is a matrix over $\mathcal{O}_{\mathbb{K}}[[\G]]$. Similarly we can find integral matrices $\mathcal{E}_\C$ and $\mathcal{E}_\D$ such that $A_0^{-1}A\frac{d\mathcal{C}}{d\G}+A_0^{-1}Y\mathcal{C}=p^\varepsilon \mathcal{E}_\C$ and $\frac{d\mathcal{D}}{d\G}BB_0^{-1}+\mathcal{D}X B_0^{-1}=p^\varepsilon \mathcal{E}_\D$. Our goal is proving that the polynomial \begin{equation*}\label{eq:errorProp}\text{$p^{-\varepsilon}(\mathcal{K}-K)\bmod \G^\ell$ has $(5\g+1,\psi)$-log convergence.}\end{equation*}
This would imply that the valuation of $(\mathcal{K}-K)\bmod\G^\ell$ is at least $-(5\g+1)\lceil\log_p\ell\rceil-\psi+\varepsilon=m$, and hence the computed $\mathcal{K}$ agrees with the actual solution $K$ modulo $(p^m,\G^\ell)$. Moreover, as $K\bmod p^m$ has degree less than $\ell$, we even have $\K\equiv K\bmod p^m$, as required. It is easy to see that this also implies that $\L_{\ell-1}\equiv K(1)\bmod p^m$.\\

We follow the proofs of Lemma 17 and 19 of \cite{HubrechtsHECEven}. Let $L$ be a matrix such that $p^\varepsilon LD=\D-D$. This gives
\[p^\varepsilon\mathcal{E}_\D=p^\varepsilon\frac{d(LD)}{d\G}BB_0^{-1}+(\D-D)XB_0^{-1} =p^\varepsilon\frac{dL}{d\G}DBB_0^{-1},\]
or $dL/d\G=\mathcal{E}_\D B_0 B^{-1}D^{-1}$. We note that $L_0=0$ and integrate to find
\[p^{-\varepsilon}(\D-D)=LD=\left(\int \mathcal{E}_\D B_0B^{-1}D^{-1}d\G\right)D.\]
By Lemma \ref{lem:convergenceLogProduct} we see that $\mathcal{E}_\D B_0B^{-1}D^{-1}$ has $(\gamma,2\alpha+\delta)$-log convergence, and as integrating is not worse than adding 1 to the logarithmic factor, we find that $p^{-\varepsilon}(\D-D)$ has $(2\gamma+1,2(\alpha+\delta)$-log convergence. Working similarly we find the same for $p^{-\varepsilon}(\mathcal{C}-C)$. Note that this implies that $\C$ and $\D$, which are polynomials of degree less than $\ell$, are both $(\gamma,\delta)$-log convergent.

From the calculation
\[\Delta(CK_0D)=A\frac{d C}{d\G}K_0DB+ACK_0\frac{dD}{d\G}B+ACK_0DX+YCK_0DB=0\]
we conclude that $K=C K_0D$.  Choose $L'$ such that $p^\varepsilon CL'K_0D=\K-\C \K_0\D$, then
\[p^{-\varepsilon}(\Delta'\K-\Delta'(\C \K_0\D))=A_0^{-1}AC\frac{dL'}{d\G}K_0DBB_0^{-1}.\]
Note that again $L_0'=0$, so that if we isolate $\frac{dL'}{d\G}$ and integrate, we find
\[L'=p^{-\varepsilon}\int C^{-1}A^{-1}A_0(\Delta'\K-\Delta'(\C \K_0\D))B_0B^{-1}D^{-1}K_0^{-1}d\G.\]
We know that $\Delta'\K=p^\varepsilon\mathcal{E}_\K$ and verify that $$\Delta'(\C \K_0\D)=p^\varepsilon(\mathcal{E}_\C \K_0\D BB_0^{-1}+A_0^{-1}A\C \K_0\mathcal{E}_\D).$$ This gives that $p^{-\varepsilon}(\K-\C \K_0\D)$ equals
\[C\left[\int C^{-1}A^{-1}\left(A_0\mathcal{E}_KB_0-A_0\mathcal{E}_C \K_0\mathcal{D}B-A\mathcal{C}\K_0\mathcal{E}_DB_0\right) B^{-1}D^{-1}d\G\right]D,\]
and hence has $(5\gamma+1,5(\alpha+\delta))$-log convergence. We conclude from
\begin{equation*}
p^{-\varepsilon}(\mathcal{K}-K)=p^{-\varepsilon} (\mathcal{K}-\mathcal{C}\K_0\mathcal{D})+p^{-\varepsilon}(\mathcal{C}-C) \K_0\mathcal{D}+ p^{-\varepsilon}C\K_0(\mathcal{D}-D)\end{equation*}
that $p^{-\varepsilon}(\mathcal{K}-K)$ has $(5\g+1,5(\alpha+\delta))=(5\g+1,\psi)$-log convergence.\\

In order to prove the theorem we need to bound the time and memory requirements of the algorithm. We will assume fast arithmetic, see e.g.\ \cite{BernsteinFastMultiplication}, which means that all basic ring operations in $\mathbb{K}$ can be performed in time essentially linear and memory linear in the object size. All elements of $\mathbb{K}$ that appear in the algorithm have valuation no less than $-\varepsilon$, hence modulo $p^\varepsilon$ all elements have bit size $\O(n\varepsilon\log_2p)$. Computing with these numbers requires then $\Ot(n\varepsilon\log_2p)$ bit operations.

If we use (\ref{eq:inductiveComp}) literally, each computation of some $\mathcal{K}_k$ requires $\O(\zeta^2)$ matrix multiplications. However, the right hand side of (\ref{eq:inductiveComp}) is essentially the coefficient of $\G^k$ in a sum of three products of (matrix) polynomials of degree $\mathcal{O}(\zeta)$, and --- using fast multiplication methods for polynomials over arbitrary algebras, see \cite{Kaminski} or \cite{CantorKaltofen} --- can thus be computed using only $\Ot(\zeta)$ matrix multiplications. As we need $K(\G)\bmod \G^\ell$ this gives in total a time requirement of $\Ot(\ell\zeta d^\omega n\varepsilon\log_2p)$ for both algorithms. Moreover, the size of $K(\G)$ determines the memory requirements for the first algorithm, hence we need $\O(\ell d^2n\varepsilon\log_2p)$ bits of space. For the second algorithm only $\O(\zeta)$ matrices over $\mathbb{K}$ have to be kept in memory, and this gives $\O(\zeta d^2n\varepsilon\log_2p)$ space.\hfill$\blacksquare$
\begin{note}\label{note:note1}
Let $\g'\in\mathcal{O}_{\mathbb{K}}$, then it is in a similar way possible to compute $K(\g')$ modulo $p^m$ with the same time and space requirements as for $K(1)$ in the theorem.
\end{note}
\begin{corollary}\label{cor:lotsOfCurves}
Let $\g_1,\ldots,\g_\ell\in\mathcal{O}_{\mathbb{K}}$ be given with accuracy $m' := m + \alpha$, then with the same assumptions as in Theorem \ref{thm:solvingDiffEq} we can compute all matrices $K(\g_1),\ldots,K(\g_\ell)$ mod $p^m$ in $\Ot(\ell \zeta d^\omega n\varepsilon\log_2p)$ bit operations and $\O(\ell  d^2n\varepsilon\log_2p)$ bits of memory.
\end{corollary}
\textsc{Proof.} We will use \emph{fast multipoint evaluation}. Let $f(x)$ be a polynomial of degree less than $\ell$ over a ring $R$. In Section 10.1 of \cite{ModernCompAlg} is explained how to evaluate $f(x)$ in $\ell$ elements of $R$ at once in such a way that it requires only $\Ot(\ell)$ arithmetic operations in $R$. Hence, taking $R=\mathcal{O}_{\mathbb{K}}$ we can compute $f(\g_1),\ldots,f(\g_\ell)$ modulo $p^{m'}$ in time $\Ot(\ell nm'\log_2p)$ and space $\O(\ell nm'\log_2p)$.

We use Theorem \ref{thm:solvingDiffEq} to compute $K(\G)$ modulo $p^{m}$. As $K(\G)$ need not be integral, we work with $p^\alpha K(\G)\bmod p^{m'}$, a matrix polynomial over $\mathcal{O}_{\mathbb{K}}$ of degree less than $\ell$. Now we can use the above result to find $p^\alpha K(\g_1), \ldots, p^\alpha K(\g_\ell)$ modulo $p^{m'}$ in time $\Ot(d^2\ell nm'\log_2p)$ and space $\O(d^2\ell n m'\log_2p)$. Taking the maximum of this result (note that $m'\leq \varepsilon$) and the complexities of Theorem \ref{thm:solvingDiffEq} concludes the proof.\ep

\section{Hyperelliptic curves in odd characteristic}\label{sec:oddChar}
In this section we will use some results from our paper \cite{HubrechtsHECOdd} about the application of deformation in point counting. Let $p$ be an odd prime and suppose we are given a hyperelliptic curve $\bar E_1$ over $\F_{p^n}$ of genus $g$ in Weierstrass form \[y^2=\bar Q_1(x) = x^{2g+1}+\sum_{i=0}^{2g}\bar a_i x^i\quad\in\F_{p^n}[x],\]
where $\bar Q_1$ is squarefree. The purpose of this section is to compute the zeta function of this curve in a memory efficient way, using Theorem \ref{thm:solvingDiffEq}. The basic idea is to deform this equation to one defined over $\Fp$, which will give us a differential equation of the kind considered in the previous section. In \cite{HubrechtsHECOdd} this was done by taking a family $y^2=\bar Q(x,\G)$ over $\Fp$ or a small extension field and substituting some $\og\in\F_{p^n}$ for $\G$. This method however does not allow us to compute the zeta function of a general hyperelliptic curve over $\F_{p^n}$. In this paper we let $\bar Q(x,\G)$ be defined over $\F_{p^n}$ and we then specialize to $\G=1$. Combining this with Theorem \ref{thm:solvingDiffEq} yields our memory efficient algorithm. We assume $p$ to be fixed in all complexity estimates of this section.

\subsection{Overview of the deformation theory}\label{ssec:overviewDeform}
Let $\bar Q_0(x) = x^{2g+1}+\sum_{i=0}^{2g}\bar b_i x^i\in\F_p[x]$ define a hyperelliptic curve $\bar E_0$ of genus $g$, for example $\bar Q_0(x):=x^{2g+1}+1$ if $p\nmid 2g+1$ and $\bar Q_0(x):=x^{2g+1}+x$ otherwise. We write $\Q_{p^n}$ for the unique unramified degree $n$ extension of $\Qp$, $\sigma$ denotes the $p$th power Frobenius automorphism on $\Q_{p^n}$ and $\Z_{p^n}$ is the ring of integers of $\Q_{p^n}$. We recall that the \emph{Teichm\"{u}ller lift} $\g\in\Z_{p^n}$ of $\bar\g\in\F_{p^n}$ is the unique root of unity that reduces to $\bar\g$ modulo $p$. Further on we will also need to extend $\sigma$ with $\sigma(\G):=\G^p$; the projection $\Z_{p^n}\to\F_{p^n}$ is always denoted with $\bar\ $ and an algebraic closure of a field $k$ is denoted as $k^\text{alg\,cl}$.

Let $a_i\in\Z_{p^n}$ and $b_i\in\Zp$ be (arbitrary) lifts of the coefficients $\bar a_i$ and $\bar b_i$, which gives us also lifts $Q_1$ and $Q_0$ of $\bar Q_1$ resp.\ $\bar Q_0$ --- monic polynomials of degree $2g+1$. We define the polynomial \[Q(x,\G) := x^{2g+1}+\sum_{i=0}^{2g}\left[(a_i-b_i)\G+b_i\right]x^i,\]
which gives a hyperelliptic curve $\bar E_{\og}\leftrightarrow y^2=\bar Q(x,\og)$ for almost all $\og\in\F_{p^n}^\text{alg\,cl}$, and makes our notation consistent regarding $\bar E_1$ and $\bar E_0$:
\[\bar E_1\  \longleftrightarrow\ y^2=\bar Q(x,1)=\bar Q_1(x)\ \ \ \text{ and }\ \ \ \bar E_0 \longleftrightarrow\ y^2=\bar Q(x,0)=\bar Q_0(x).\]
We now give a short overview of the theory in \cite{HubrechtsHECOdd}. Let $r(\G)$ be the resultant
\[r(\G) := \text{Res}_x\left(Q(x,\G); \frac{\partial}{\partial x}Q(x,\G)\right),\]
then it is clear from the construction of $Q(x,\G)$ that $r(0)$ and $r(1)$ are units in $\Z_{p^n}$ and $\rho := \deg r(\G)\leq 4g$. Suppose $r(\G)=\sum_{i=0}^\rho r_i\G^i$; with $\rho'$ the degree of $\bar r(\G)$ we define $\tilde r(\G):=\sum_{i=0}^{\rho'}r_i\G^i$. We defined in Sections 3.2 and 3.3 of \cite{HubrechtsHECOdd} a ring $S$ and an $S$-module $T$, which can be represented as (where $\dagger$ means overconvergent completion):
\begin{align*}
S &:= \Q_{p^n}\left[\G,\tilde r(\G)^{-1}\right]^\dagger,\\
T &:= \cfrac{\Q_{p^n}\left[x,y,y^{-1},\G,\tilde r(\G)^{-1}\right]^\dagger} {(y^2-Q(x,\G))}.\end{align*}
Let $d:T\to Tdx$ be the differential $\frac{\partial}{\partial x}dx$ and $\nabla:T\to Td\G$ the connection $\frac{\partial}{\partial \G}d\G$ such that $d(\G)=\nabla(x)=0$. Then we showed that a certain submodule $H_{MW}^-$ of $Tdx/d(T)$ is a free $S$-module of rank $2g$, which, after substituting for $\G$ any Teichm\"{u}ller lift $\g\in\Q_{p^n}$ which is no zero modulo $p$ of $r(\G)$, gives the same $2g$-dimensional $\Q_{p^n}$-vector space as Kedlaya's $A^\dagger\otimes\Q_{p^n}$, defined in Section 3 of \cite{KedlayaCountingPoints}. We also constructed a Frobenius map $F_p$ on $H_{MW}^-$ which after the specialization $\G\leftarrow\g$ again equals Kedlaya's. The following diagram is well defined and commutes:
\begin{equation}\label{eq:diagram}\begin{CD}
H_{MW}^- @>{\nabla}>> H_{MW}^-d\G\\
@VV{F_p}V @VV{F_p}V\\
H_{MW}^- @>{\nabla}>> H_{MW}^-d\G.
\end{CD}\end{equation}
We have the $S$-basis $\{\frac{x^idx}{\sqrt{Q}}\}_{i=0}^{2g-1}$ for $H_{MW}^-$ and can hence define $(2g\times 2g)$-matrices over $S$ for our operators, namely $G(\G)$ for $\nabla$ and $F(\G)$ for $F_p$, e.g.\ $F_p(x^idx/\sqrt{Q})=\sum_j F_{ij}(\G)x^jdx/\sqrt{Q}$. As Kedlaya showed in \cite{KedlayaCountingPoints}, the main step in computing the zeta function of $\bar E_1$ is to compute $F(1)$ up to a certain precision.

Let $H(\G):=r(\G)G(\G)$, then the equation, derived from (\ref{eq:diagram}), at the end of Section 3.6 in \cite{HubrechtsHECOdd} reads
\[rr^\sigma\frac{dF}{d\G}+r^\sigma FH-p\G^{p-1}rH^\sigma F=0.\]
Here we use $r^\sigma$ and $H^\sigma$ for $r^\sigma(\G^p)$ respectively $H^\sigma(\G^p)$.
Substituting $K:=r^MF$ for some integer $M\geq 0$ that is made more precise below, this becomes
\[r^\sigma \frac{dK}{d\G} r + r^\sigma K(H-M\frac{d r}{d\G})+(-p\G^{p-1}H^\sigma)Kr=0,\]
which is of the form $\Delta K=0$ explained in Section \ref{ssec:diffEq} above, with $d=2g$, $A=r^\sigma$, $B=r$, $X=H-M\frac{d r}{d\G}$ and $Y=-p\G^{p-1}H^\sigma$.

\subsection{Computing the zeta function}\label{ssec:computeZeta}
We will now determine the constants $\ell,\zeta,m$ and $\varepsilon$ in order to apply Theorem \ref{thm:solvingDiffEq}. From Proposition 16 and Lemma 18 in \cite{HubrechtsHECOdd} it follows that with $\alpha:=(2g-1)(\log_p g +2)+g = \O(g\log g)$ we have
\[\ord(F)=\ord(K)\geq -\alpha\quad\text{ and }\quad \ord(F^{-1})=\ord(K^{-1})\geq -\alpha.\]
Proposition 17 of \cite{HubrechtsHECOdd} (with $\kappa:=\deg_{\G}Q(x,\G)=1$) shows that $\deg H\leq 8g$ and as a consequence we can take
\[\zeta := \max\{(p+1)\rho,p\rho+8g+1,p+8pg+\rho\} = \O(g).\]
We note in passing that this Proposition 17 also implies that $\ord(H)\geq\frac{-10g}{p-1}$ and hence the conditions $\ord(X),\ord(Y)\geq -\psi$ at the end of Section \ref{ssec:assumptions} will be met.

We need $F(1)$ modulo $p^m$ with $m$ defined as $N_b$ in Section 4 of \cite{HubrechtsHECOdd}, namely (with $a=1$)
\[m := \left\lceil \frac{ng}2+(2g+1)\log_p2\right\rceil+ n\left\lfloor\log_p(g)+2\right\rfloor+\lfloor 2gn(\log_pg+3)\rfloor=\O(ng\log g).\]
The exponent $M:=p(2m+4)+(p-1)/2=\O(ng\log g)$ of $r(\G)$ and the precision $\ell$ are given by Proposition 16 in \cite{HubrechtsHECOdd}:
\[\ell := (2m+5)(8g+2)p+1=\O(ng^2\log g).\]
Next we need $\g$ and $\delta$ such that (\ref{eq:CDlogConv}) holds. For the solution $C$ of $A\frac{dC}{d\G}+YC=0$ we can find this in Proposition 20 of \cite{HubrechtsHECOdd}: the matrix $C$ in that proposition does not correspond to $C$ in this paper, but the result and proof are completely the same. The conclusion is that
\[C\text{ and }C^{-1}\text{ have }(2g\log_pg+g,0)\text{-log  convergence.}\]
We have that $K=CK_0D$, and as a consequence $D=K_0^{-1}C^{-1}K$ and $D^{-1}=K^{-1}CK_0$ have $(2g\log_pg+g,2\alpha)$-log convergence. Hence we can take $\g:=2g\log_pg+g=\O(g\log g)$ and $\delta:=2\alpha$. Now with $\psi = 2\alpha+5\delta=12\alpha$ we find
\[\varepsilon = m+(5\g+1)\lceil\log_p\ell\rceil+\psi = \O(ng(\log g)^2).\]
The analysis in \cite{HubrechtsHECOdd}, namely Steps 1, 2 and 5 of Section 6.3, shows that the time and space requirements for computing $r$, $H$ and $K_0$ will not have any influence on the result, and as a consequence we can apply Theorem \ref{thm:solvingDiffEq} to find $K(1)\bmod p^m$ in time
\[\Ot(\ell\zeta g^\omega n\varepsilon) = \Ot(g^{4+\omega}n^3)\]
and with memory requirements
\begin{equation}\label{eq:memoryOdd}
\O(\zeta g^2n\varepsilon)=\O(g^4(\log g)^2n^2).\end{equation}

To conclude the algorithm we still need to approximate the matrix $\mathcal{F}$ of $F_p^n$. First we compute $F(1)=r(1)^{-M}K(1)$ and then
\[\mathcal{F} = F(1)^{\sigma^{n-1}}\cdot F(1)^{\sigma^{n-2}}\cdots F(1)^\sigma\cdot F(1),\]
which can be certainly done in time $\Ot(g^3n^3+g^{1+\omega}n^2)$ and memory $\O(n^2g^3)$ as explained in \cite{KedlayaCountingPoints} (see however Section \ref{ssec:manyCurves} for a much faster method). The numerator of the zeta function equals $\det(1-\mathcal{F}t)$ and can be found in time $\Ot(g^{2+\omega}n^2)$. These last complexities can all be found in Step 8 of Section 6.3 of \cite{HubrechtsHECOdd}, where the memory requirements are bounded by (\ref{eq:memoryOdd}). This results in the following theorem.
\begin{theorem}\label{thm:OddChar}
There exists an explicit and deterministic algorithm to compute the zeta function of any hyperelliptic curve of genus $g$ over $\F_{p^n}$, with $p$ odd, that uses $\Ot(n^3g^{4+\omega})$ bit operations and bit space $\O(n^2g^4(\log_pg)^2)$.\end{theorem}

\section{Additional applications}\label{sec:manyCurves1}
\subsection{Hyperelliptic curves in even characteristic}\label{sec:evenChar}

By an argument similar to the one explained in the previous section, we can prove the following result.
\begin{theorem}\label{thm:EvenChar}
There exists an explicit and deterministic algorithm that computes the zeta function of any hyperelliptic curve of genus $g$ over $\F_{2^n}$ using  $\Ot(n^3g^{4+\omega+3\tau})$ bit operations and $\O(n^2g^{4+\tau}(\log g)^{2+\tau})$ bits of memory. Here $\tau=0$ for almost all curves and $\tau=1$ in the general case.\end{theorem}
In order to show this, one uses the results from \cite{HubrechtsHECEven}, in turn partly inspired by Denef and Vercauteren's article \cite{DenefVercauteren}. We want to point out that $\tau$ in the theorem above will be 0 precisely when all $m_i$ are equal to 1 (or are bounded by $\mathcal{O}(1)$) in the notation of the beginning of Section 4 of \cite{DenefVercauteren}. More details can be found in Section 5.3 of \cite{HubrechtsThesis}.

\subsection{Many curves at once}\label{ssec:manyCurves}

If we choose a family defined over $\F_{p^n}$ as in Section \ref{ssec:overviewDeform} (or in an analogous way in characteristic 2), e.g.\ given by $\bar E_{\G}:y^2=\bar Q(x,\G)$, and $\og_1,\ldots,\og_{g^2n}\in\F_{p^n}$, we can compute the zeta functions of the curves $\bar E_{\og_i}$ all at once in a very efficient way. There are three main steps needed in order to achieve this. First, computing the Teichm\"{u}ller lifts of all $\bar \g_i$ modulo $p^\varepsilon$ can be done in time $\O(g^2n(n\varepsilon)^{1+\rho}) or \O((n^3g^3)^{1+\rho})$ for any $\rho>0$ as shown in
\cite[Proposition 6]{HubrechtspAdicArithmetic}. Second, we compute all the matrices of the $p$th power Frobenius using Corollary \ref{cor:lotsOfCurves} above in time $\Ot(n^3g^{5+\omega})$. And third, in order to retrieve the matrices of the $q$th power Frobenius and hence the zeta functions, we can use Kedlaya's trick \cite[Section 5]{KedlayaCountingPoints} combined with Proposition 3 of \cite{HubrechtspAdicArithmetic}, resulting in a time complexity bounded by $\O((g^2n)(n^2g^{1+\omega})^{1+\rho})$. Noting that all these algorithms are deterministic we conclude with the following theorem:
\begin{theorem}\label{thm:ManyCurves}Suppose we are given a family $\bar E_\G:y^2=\bar Q(x,\G)$ over $\F_{p^n}$ and $\bar\g_1,\ldots,\bar\g_{g^2n}\in\F_{p^n}$ such that $\bar E_0$ is defined over $\F_p$ and all $\bar E_{\bar \g_i}$ and $\bar E_0$ are hyperelliptic curves of genus $g$. Choose $\rho>0$. There exists an explicit and deterministic algorithm that computes the zeta functions of all curves $\bar E_{\bar \g_i}$ that requires $\Ot(n^{3+\rho}g^{5+\omega})$ bit operations.\end{theorem}
In an obvious way a similar algorithm can be shown to exist for characteristic 2.

This result could be interesting if one wants to find a curve with a special property, as is the case in cryptography. For example, suppose that we want to find a curve over $\F_{p^n}$ with $N$ as order of its jacobian, such that $N$ has a very large prime factor. Then we can expect that we have to try $\O(n)$ curves in order to find such a curve, and this is exactly something we can do very efficiently as explained above.\\

In \cite{CastryckHubrechtsVercauteren} deformation is used for the computation of the zeta function of $C_{a,b}$ curves, and as explained in Section 5.4 of that paper similar results as above apply.

\subsection{Hypersurfaces}\label{ssec:hypersurfaces}
Finally we can also use our memory efficient algorithm for solving differential equations in the context of hypersurfaces. For example, Lauder gives in \cite{LauderDeformation} an algorithm that computes the zeta function of certain hypersurfaces satisfying an `almost diagonal' equation over $\F_{p^n}$. As he uses deformation as the main step in this result, the memory requirements drop from cubic to quadratic in $n$ using the result in this paper. Gerkmann discusses in \cite{GerkmannHypersurfaces} several deformation strategies for smooth projective surfaces, and an important step in there is again solving such a differential equation. Although the improvements depend on the type of algorithm considered, most algorithms presented in \cite{GerkmannHypersurfaces} will profit from Theorems \ref{thm:solvingDiffEq} and \ref{thm:ManyCurves}.

\bibliographystyle{amsplain}
\bibliography{bibliography}

\providecommand{\bysame}{\leavevmode\hbox to3em{\hrulefill}\thinspace}
\providecommand{\MR}{\relax\ifhmode\unskip\space\fi MR }
\providecommand{\MRhref}[2]{%
  \href{http://www.ams.org/mathscinet-getitem?mr=#1}{#2}
}
\providecommand{\href}[2]{#2}
\begin{thebibliography}{10}

\bibitem{BernsteinFastMultiplication}
Daniel~J. Bernstein, \emph{Fast multiplication and its applications},
  Algorithmic number theory: lattices, number fields, curves and cryptography,
  Math. Sci. Res. Inst. Publ., vol.~44, Cambridge Univ. Press, Cambridge, 2008,
  pp.~325--384.

\bibitem{CantorKaltofen}
David~G. Cantor and Erich Kaltofen, \emph{On fast multiplication of polynomials
  over arbitrary algebras}, Acta Inform. \textbf{28} (1991), no.~7, 693--701.

\bibitem{CastryckHubrechtsVercauteren}
Wouter Castryck, Hendrik Hubrechts, and Frederik Vercauteren, \emph{Computing
  zeta functions in families of {$C\sb {a,b}$} curves using deformation},
  Algorithmic number theory, Lecture Notes in Comput. Sci., vol. 5011,
  Springer, Berlin, 2008, pp.~296--311.

\bibitem{CohenFrey}
Henri Cohen, Gerhard Frey, Roberto Avanzi, Christophe Doche, Tanja Lange, Kim
  Nguyen, and Frederik Vercauteren (eds.), \emph{Handbook of elliptic and
  hyperelliptic curve cryptography}, Discrete Mathematics and its Applications
  (Boca Raton), Chapman \& Hall/CRC, Boca Raton, FL, 2006.

\bibitem{CoppersmithWinograd}
Don Coppersmith and Shmuel Winograd, \emph{Matrix multiplication via arithmetic
  progressions}, J. Symbolic Comput. \textbf{9} (1990), no.~3, 251--280.

\bibitem{DVErratum}
Jan Denef and Frederik Vercauteren, \emph{Errata for ``{A}n extension of
  {K}edlaya's algorithm to hyperelliptic curves in characteristic 2'', and
  related papers}, Available on
  \texttt{http://www.wis.kuleuven.be/algebra/denef{\_}\allowbreak{}papers/Erra%
taPointCounting.pdf}.

\bibitem{DenefVercauteren}
\bysame, \emph{An extension of {K}edlaya's algorithm to hyperelliptic curves in
  characteristic 2}, J. Cryptology \textbf{19} (2006), no.~1, 1--25, Erratum
  available as \cite{DVErratum}.

\bibitem{GerkmannHypersurfaces}
Ralf Gerkmann, \emph{Relative rigid cohomology and deformation of
  hypersurfaces}, Int. Math. Res. Pap. IMRP (2007), no.~1, Art. ID rpm003, 67.

\bibitem{GerkmannEC}
\bysame, \emph{Relative rigid cohomology and point counting on families of
  elliptic curves}, J. Ramanujan Math. Soc. \textbf{23} (2008), no.~1, 1--31.

\bibitem{HubrechtsThesis}
Hendrik Hubrechts, \emph{Elliptic and hyperelliptic curve point counting
  through deformation}, PhD thesis, KULeuven, Belgium. May 2007. Available on
  \texttt{http://wis.kuleuven.be/algebra/hubrechts/}.

\bibitem{HubrechtspAdicArithmetic}
Hendrik Hubrechts, \emph{Fast arithmetic in unramified $p$-adic fields}, To
  appear in Finite Fields and Their Applications. Preprint available on
  {\verb!http://wis.kuleuven.be/algebra/hubrechts/!}.

\bibitem{HubrechtsHECEven}
\bysame, \emph{Point counting in families of hyperelliptic curves in
  characteristic 2}, LMS J. Comput. Math. \textbf{10} (2007), 207--234
  (electronic).

\bibitem{HubrechtsHECOdd}
\bysame, \emph{Point counting in families of hyperelliptic curves}, Found.
  Comput. Math. \textbf{8} (2008), no.~1, 137--169.

\bibitem{Kaminski}
Michael Kaminski, \emph{An algorithm for polynomial multiplication that does
  not depend on the ring constants}, J. Algorithms \textbf{9} (1988), no.~1,
  137--147.

\bibitem{KedlayaCountingPoints}
Kiran~S. Kedlaya, \emph{Counting points on hyperelliptic curves using
  {M}onsky-{W}ashnitzer cohomology}, J. Ramanujan Math. Soc. \textbf{16}
  (2001), no.~4, 323--338.

\bibitem{LauderDeformation}
Alan G.~B. Lauder, \emph{Deformation theory and the computation of zeta
  functions}, Proc. London Math. Soc. (3) \textbf{88} (2004), no.~3, 565--602.

\bibitem{TsuzukiKloosterman}
N.~Tsuzuki, \emph{Bessel {F-isocrystals} and an algorithm of computing
  {Kloosterman} sums}, Unpublished.

\bibitem{ModernCompAlg}
Joachim von~zur Gathen and J{\"u}rgen Gerhard, \emph{Modern computer algebra},
  Cambridge University Press, Cambridge, 2003.

\end{thebibliography}

\end{document}